\documentclass{llncs}

\usepackage{amsmath,amsfonts,amssymb}

\usepackage{tikz}
\usetikzlibrary{backgrounds}

\usepackage{tikzit}

\tikzstyle{thick}=[-, draw=black, tikzit draw=black, thick]

\usepackage{svg}
\usepackage{graphics}

\usepackage{hyperref}

\newcommand\id{\mathrm{id}}
\newcommand\src{\mathfrak{s}}
\newcommand\tgt{\mathfrak{t}}

\title{Rendering string diagrams recursively}
\author{Celia Rubio-Madrigal\inst{1} \and Jules Hedges\inst{2,3}}
\institute{CISPA Helmholtz Center for Information Security \and Mathematically Structured Programming group, University of Strathclyde \and Institute for Categorical Cybernetics}

\begin{document}

\maketitle

\begin{abstract}
	String diagrams are a graphical language used to represent processes that can be composed sequentially or in parallel, which correspond graphically to horizontal or vertical juxtaposition. In this paper we demonstrate how to compute the layout of a string diagram by folding over its algebraic representation in terms of sequential and parallel composition operators. The algebraic representation can be seen as a term of a free monoidal category or a proof tree for a small fragment of linear logic. This contrasts to existing non-compositional approaches that use graph layout techniques. The key innovation is storing the diagrams in binary space-partition trees, maintaining a right-trapezoidal shape for the diagram's outline as an invariant.

	We provide an implementation in Haskell, using an existing denotational graphics library called Diagrams. Our renderer also supports adding semantics to diagrams to serve as a compiler, with matrix algebra used as an example.
\end{abstract}

\section{Introduction}

String diagrams are a graphical language used to represent any process that can be composed sequentially or run in parallel, which translates graphically to joining them horizontally or vertically. Algebraically, such processes form a (symmetric) monoidal category, and string diagrams are a presentation of free monoidal categories. They have many applications, including but not limited to electrical circuits \cite{fong_algebra_2016}, probability theory \cite{fritz_synthetic_2020}, game theory \cite{ghani_compositional_2018}, machine learning \cite{foundations-gradient-learning} and abstract algebra \cite{interacting-bialgebras}. The most successful applications have been in quantum mechanics \cite{chapman_picturing_2018} and natural language processing \cite{discocat}.

String diagrams have an underlying ``algebraic'' representation with operations for sequential and parallel composition, with an equational theory referred to as a \emph{monoidal category}. This paper considers the problem of converting from the algebraic representation to a ``rendered'' string diagram. Conceptually, this is straightforward: it is a \emph{fold} over the datatype where sequential and tensor composition are recursively translated to horizontal and vertical juxtaposition. However, building an actual implementation, especially one that is able to produce aesthetically pleasing diagrams, involves more intricate geometric details. Although string diagrams have usually been considered as contained within a rectangle, as far back as the foundational work of Joyal and Street \cite{joyal-street}, our key insight is that they should be contained inside a \emph{trapezoid}. These geometric details, and an implementation in Haskell, are the contribution of this paper.

\textbf{Acknolwedgements.} This paper is based on the first author's MSc thesis \cite{celia-thesis}, which was supported by a fellowship from ''la Caixa'' Foundation (ID 100010434, with fellowship code LCF/BQ/EU22/11930080).

\section{Free monoidal categories}

A \emph{monoidal signature} \cite{joyal-street} $\mathcal S$ consists of a set $\mathrm{Ob} (\mathcal S)$ of \emph{object symbols}, a set $\mathrm{Mor} (\mathcal S)$ of \emph{morphism symbols}, and, for each morphism symbol $f\in\mathrm{Mor} (\mathcal S)$, a pair of finite ordered lists of object symbols, $\src (f)$ and $\tgt (f)$, called the \emph{source} and \emph{target} of $f$. We think of a list $\left< x_1, \ldots, x_n \right>$ as a formal tensor product $x_1 \otimes \cdots \otimes x_n$. Monoidal signatures are closely related to directed hypergraphs and Petri nets.

We are interested in representing the free monoidal category generated by a monoidal signature, $F_{MC} (\mathcal S)$. Closely related to this, we can generate other classes of structured monoidal categories using the same generating data, such as the free symmetric monoidal category $F_{SMC} (\mathcal S)$, the free traced monoidal category $F_{TMC} (\mathcal S)$ and the free compact closed category $F_{CCC} (\mathcal S)$ \cite{selinger-survey}.

Two common representations of the free monoidal category generated by a signature are \emph{terms} and \emph{string diagrams}. A third representation called \emph{brick diagrams} was introduced in \cite{foundations-brick-diagrams} that bridges the gap between the two.

\subsection{Terms}

The class of \emph{terms} is defined recursively, together with an extension of the functions $\src$ and $\tgt$ associating a list of object generators to each term:
\begin{itemize}
	\item Every morphism symbol is a term.
	\item For every list of object symbols $x$ there is a term $\id_x$ with $\src (\id_x) = \tgt (\id_x) = x$.
	\item For every pair of terms $\alpha, \beta$ with $\tgt (\alpha) = \src (\beta)$ there is a term $\alpha; \beta$ such that $\src (\alpha; \beta) = \src (\alpha)$ and $\tgt (\alpha; \beta) = \tgt (\beta)$.
	\item For every pair of terms $\alpha, \beta$ there is a term $\alpha \otimes \beta$ such that $\src (\alpha \otimes \beta) = \src (\alpha) + \src (\beta)$ and $\tgt (\alpha \otimes \beta) = \tgt (\alpha) + \tgt (\beta)$ (where $+$ is list concatenation).
\end{itemize}

We refer to the lengths of the lists $\src (\alpha)$ and $\tgt (\alpha)$ as the term's source/target \emph{arity}. Equivalently, every term $\alpha$ can be seen as a proof of the sequent $\src (\alpha) \vdash \tgt (\alpha)$ in a small fragment of noncommutative linear logic consisting of three proof rules:
\[ \frac{}{x \vdash x}\text{(Ax)} \qquad \frac{s \vdash t \qquad s' \vdash t'}{s + s' \vdash t + t'}\text{($\otimes$)} \qquad \frac{s \vdash t \qquad t \vdash u}{s \vdash u}\text{(cut)} \]
Note that this fragment does not allow cut elimination. 

The term representation has the advantage of simplicity, but the disadvantage that there are multiple terms representing the same morphism of a free monoidal category. The equations between terms we need to consider are exactly the ones generated by the following rules: (1) $;$ is associative and unital with units being $\id$ of the appropriate type, (2) $\otimes$ is associative and unital with unit $\id_\epsilon$, and (3) the \emph{interchange law} $(\alpha; \beta) \otimes (\gamma; \delta) = (\alpha \otimes \gamma); (\beta \otimes \delta)$ for all terms $\alpha, \beta, \gamma, \delta$ satisfying $\tgt (\alpha) = \src (\beta)$ and $\tgt (\gamma) = \src (\delta)$.
Of these, the unitality and associativity of $;$ and $\otimes$ can be trivialised by passing to \emph{unbiased} terms, ie. using $n$-ary rather than binary composition. However, the interchange law is fundamentally difficult.

\subsection{String diagrams}

String diagrams are a graph-like topological representation of terms. They satisfy a fundamental \emph{coherence theorem} \cite{joyal-street}: that the above equational theory for terms coincides with \emph{planar isotopy} of diagrams; in particular, they trivialise the interchange law. We will illustrate them by the following example.

\begin{figure}[t]
\[ \begin{tikzpicture}[scale=0.7]
	\begin{pgfonlayer}{nodelayer}
		\node [style=none] (0) at (-4, 4) {};
		\node [style=none] (1) at (4, 4) {};
		\node [style=none] (2) at (-4, -4) {};
		\node [style=none] (3) at (4, -4) {};
		\node [style=none] (4) at (-2, 2) {};
		\node [style=none] (5) at (2, 1) {};
		\node [style=none] (6) at (2, -2) {};
		\node [style=none] (7) at (-2, -1) {};
		\node [style=none] (8) at (-4, 2) {};
		\node [style=none] (9) at (-4, -2) {};
		\node [style=none] (10) at (4, 0) {};
		\node [style=none] (11) at (4, 2) {};
		\node [style=none] (12) at (4, -2) {};
		\node [style=none] (13) at (-4, 0) {};
		\node [style=none] (14) at (-2, 2.75) {$f$};
		\node [style=none] (15) at (2, 1.75) {$f$};
		\node [style=none] (16) at (-2, -0.25) {$g$};
		\node [style=none] (17) at (2, -1.25) {$g$};
	\end{pgfonlayer}
	\begin{pgfonlayer}{edgelayer}
		\draw (0.center) to (1.center);
		\draw (1.center) to (3.center);
		\draw (3.center) to (2.center);
		\draw (2.center) to (0.center);
		\draw (8.center) to (4.center);
		\draw [in=135, out=0, looseness=1.25] (13.center) to (7.center);
		\draw [in=-135, out=0, looseness=1.25] (9.center) to (7.center);
		\draw [in=180, out=30, looseness=1.25] (4.center) to (5.center);
		\draw [in=135, out=-45] (4.center) to (6.center);
		\draw [in=-135, out=0, looseness=1.25] (7.center) to (6.center);
		\draw [in=180, out=45, looseness=1.50] (5.center) to (11.center);
		\draw [in=-180, out=-45, looseness=1.25] (5.center) to (10.center);
		\draw (6.center) to (12.center);
	\end{pgfonlayer}
\end{tikzpicture} \qquad\qquad \begin{tikzpicture}[scale=0.7]
	\begin{pgfonlayer}{nodelayer}
		\node [style=none] (0) at (-4, 4) {};
		\node [style=none] (1) at (4, 4) {};
		\node [style=none] (2) at (-4, -4) {};
		\node [style=none] (3) at (4, -4) {};
		\node [style=none] (8) at (-4, 2) {};
		\node [style=none] (9) at (-4, -2) {};
		\node [style=none] (10) at (4, 0) {};
		\node [style=none] (11) at (4, 2) {};
		\node [style=none] (12) at (4, -2) {};
		\node [style=none] (13) at (-4, 0) {};
		\node [style=none] (14) at (-2, 2) {$f$};
		\node [style=none] (15) at (-2, -1) {$g$};
		\node [style=none] (16) at (2, 1) {$f$};
		\node [style=none] (17) at (2, -2) {$g$};
		\node [style=none] (18) at (-2.5, 2.75) {};
		\node [style=none] (19) at (-1.5, 2.75) {};
		\node [style=none] (20) at (-1.5, 1.25) {};
		\node [style=none] (21) at (-2.5, 1.25) {};
		\node [style=none] (22) at (1.5, 1.75) {};
		\node [style=none] (23) at (2.5, 1.75) {};
		\node [style=none] (24) at (2.5, 0.25) {};
		\node [style=none] (25) at (1.5, 0.25) {};
		\node [style=none] (26) at (-2.5, -0.25) {};
		\node [style=none] (27) at (-1.5, -0.25) {};
		\node [style=none] (28) at (-1.5, -1.75) {};
		\node [style=none] (29) at (-2.5, -1.75) {};
		\node [style=none] (30) at (1.5, -1.25) {};
		\node [style=none] (31) at (2.5, -1.25) {};
		\node [style=none] (32) at (2.5, -2.75) {};
		\node [style=none] (33) at (1.5, -2.75) {};
		\node [style=none] (34) at (-2.5, 2) {};
		\node [style=none] (35) at (1.5, 1) {};
		\node [style=none] (36) at (-1.5, -1) {};
		\node [style=none] (37) at (2.5, -2) {};
		\node [style=none] (38) at (-1.5, 2.25) {};
		\node [style=none] (39) at (-1.5, 1.75) {};
		\node [style=none] (40) at (2.5, 1.25) {};
		\node [style=none] (41) at (2.5, 0.75) {};
		\node [style=none] (42) at (-2.5, -0.75) {};
		\node [style=none] (43) at (-2.5, -1.25) {};
		\node [style=none] (44) at (1.5, -1.75) {};
		\node [style=none] (45) at (1.5, -2.25) {};
	\end{pgfonlayer}
	\begin{pgfonlayer}{edgelayer}
		\draw (0.center) to (1.center);
		\draw (1.center) to (3.center);
		\draw (3.center) to (2.center);
		\draw (2.center) to (0.center);
		\draw (18.center) to (19.center);
		\draw (19.center) to (20.center);
		\draw (20.center) to (21.center);
		\draw (21.center) to (18.center);
		\draw (22.center) to (23.center);
		\draw (23.center) to (24.center);
		\draw (24.center) to (25.center);
		\draw (25.center) to (22.center);
		\draw (26.center) to (27.center);
		\draw (27.center) to (28.center);
		\draw (28.center) to (29.center);
		\draw (29.center) to (26.center);
		\draw (30.center) to (31.center);
		\draw (31.center) to (32.center);
		\draw (32.center) to (33.center);
		\draw (33.center) to (30.center);
		\draw (8.center) to (34.center);
		\draw (37.center) to (12.center);
		\draw [in=180, out=0] (38.center) to (35.center);
		\draw [in=180, out=0, looseness=0.75] (39.center) to (44.center);
		\draw [in=-180, out=0] (40.center) to (11.center);
		\draw [in=180, out=0] (41.center) to (10.center);
		\draw [in=-180, out=0] (13.center) to (42.center);
		\draw [in=-180, out=0, looseness=1.25] (9.center) to (43.center);
		\draw [in=180, out=0] (36.center) to (45.center);
	\end{pgfonlayer}
\end{tikzpicture} \]
\caption{String diagram for $(f \otimes g); (f \otimes g)$ as a topological graph (left) and with standard notation (right)}
\label{fig:example-diagram}
\end{figure}
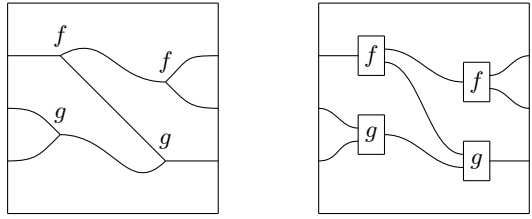

Suppose we have a monoidal signature with a single object symbol $x$, and morphism symbols $f : x \to x \otimes x$ and $g : x \otimes x \to x$. The tensor product of these is $f \otimes g : x \otimes x \otimes x \to x \otimes x \otimes x$, so we can consider the term $(f \otimes g); (f \otimes g) : x \otimes x \otimes x \to x \otimes x \otimes x$. The string diagram corresponding to this last term is depicted in \autoref{fig:example-diagram}. String diagrams are formally \emph{topological graphs}; that is, graphs equipped with a choice of planar isotopy class of topological embeddings in the plane \cite{joyal-street}, satisyfing some additional conditions. However, nodes of the graph (which are labelled by morphism symbols) are not typically depicted as point-like, but as extended regions such as rectangles for ease of readability.

\subsection{Brick diagrams}

\begin{figure}[t]
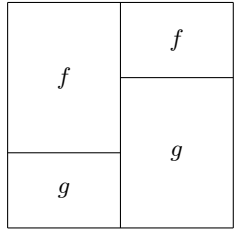

\ctikzfig{diagrams/example-brick}
\caption{Brick diagram for $(f \otimes g); (f \otimes g)$}
\label{fig:example-brick}
\end{figure}

\emph{Brick diagrams} are an alternative to string diagrams introduced by the Statebox team for use in a now-defunct diagram editor, and formally defined in \cite{foundations-brick-diagrams}.
Brick diagrams are a certain Poincar\'e dual of string diagrams, with morphism symbols represented by rectangular regions of the plane and connecting object symbols represented by lines on which they overlap. For example, the term $(f \otimes g); (f \otimes g)$ from the previous section is depicted as a brick diagram in \autoref{fig:example-brick}.
Foundationally, brick diagrams can be seen as \emph{tiling diagrams} \cite{dawson-pare} after viewing monoidal categories as degenerate \emph{double categories}.

In the conclusion of \cite{foundations-brick-diagrams} it was suggested that rendering string diagrams can be reduced to rendering brick diagrams if each individual brick is responsible for rendering a string diagram containing a single node inside itself. This idea was the starting point of this paper.

\subsection{Naive translation from terms to diagrams}

It is conceptually straightforward to convert from terms to string diagrams in a recursive way. Every generating morphism is translated to a string diagram containing a single node. The $;$ and $\otimes$ operators translate to joining side-by-side the two string diagrams that correspond to the two respective subterms, either in the horizontal or the vertical axis.
The most obvious approach is to arrange each string diagram into a square, for example $[0,1]^2$.
For compositions, we can place the sub-diagrams side by side and then use a linear transformation to squash the result back into a square.

Clearly this method can produce un-aesthetic results, with some parts of the diagram taking up far more space than others depending on the composition depth of the corresponding terms. For example, for a nested composition $f_1; (f_2; (f_3; \cdots f_n) \cdots)$, the nodes will get exponentially closer together in geometric sequence. This could be partially overcome with unbiased composition.

\begin{figure}[t]
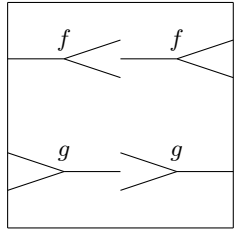

\ctikzfig{diagrams/naive-broken}
\caption{Result of applying naive composition to $(f \otimes g); (f \otimes g)$}
\label{fig:broken-diagram}
\end{figure}

More subtly, the naive method fails to produce diagrams that are even topologically correct, because we fail to preserve the intended invariant that strings are spaced equidistantly along the boundaries, and so a $;$ composition can result in a misaligned boundary.
For example, applying the above rules mechanically to the example term $(f \otimes g); (f \otimes g)$ from the previous section results in a mis-formed diagram, depicted in \autoref{fig:broken-diagram}.

\section{Trapezoid-shaped diagrams}

The key idea of this paper is that a diagram should be contained inside a \emph{right trapezoid} with side lengths determined by its source arity $\ell$ and target arity $r$. In particular, we maintain the following invariants for all diagrams:
\begin{enumerate}
    \item Diagrams must have the shape of a right trapezoid, with both right angles at the bottom.
    \item The lengths of the trapezoid's left and right sides are $\ell$ and $r$
    \item The coordinate origin is at the trapezoid's bottom left vertex.
\end{enumerate}
Note that the width of diagrams is not constrained by these invariants.

For the purposes of this paper, we are going to require that all diagrams have left and right arity $\geq 1$. Breaking one of these conditions causes the trapezoid to degenerate into a triangle, and breaking both causes it to collapse into a line segment, destroying all visual information inside. Fixing this restriction is left for future work.

\subsection{Sequential composition}

\begin{figure}[t]
    \centering
    \begin{tikzpicture}[scale=0.8]
        \draw (0,0) -- (0,1) -- (1,3) -- (1,0) -- cycle;
        \draw (1,0) -- (1,3) -- (2,2) -- (2,0) -- cycle;
        \draw[dashed] (0,1) -- (2,2);
        \draw[->] (2.5,1) to[out=60, in=120] node[midway,above] {$Pinch_r$} (3.5,1);
    \begin{scope}[shift={(4,0)}]
        \draw (0,0) -- (0,1) -- (1,1.5) -- (1,0) -- cycle;
        \draw (1,0) -- (1,3) -- (2,2) -- (2,0) -- cycle;
        \draw[dashed] (1,1.5) -- (2,2);
        \draw[->] (2.5,1) to[out=60, in=120] node[midway,above] {$Pinch_\ell$} (3.5,1);
    \end{scope}
    \begin{scope}[shift={(8,0)}]
        \draw (0,0) -- (0,1) -- (1,1.5) -- (1,0) -- cycle;
        \draw (1,0) -- (1,1.5) -- (2,2) -- (2,0) -- cycle;
    \end{scope}
    \end{tikzpicture}
    \caption{Pinch operations needed to compose two right trapezoids}
    \label{fig:composing-brickdiagrams}
    \end{figure}
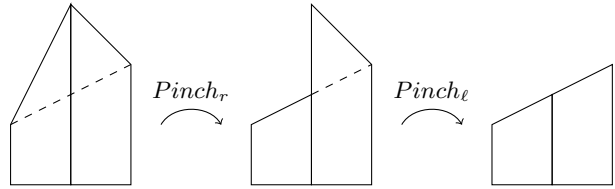

When horizontally composing two diagrams (D1;D2), the right arity of D1 equals the left arity of D2. A side-by-side placement creates matching internal connections between them, with the $i$th join occurring at point $(w_1, \frac{1}{2} + i)$.
However, the resulting shape is generally a pentagon rather than a quadrilateral. To obtain a trapezoid we use a ``pinching'' transformation that moves one top vertex of a trapezoid to a different height. As illustrated in \autoref{fig:composing-brickdiagrams}, pinching both the right corner of the left diagram and the left corner of the right diagram is required to reshape the composition back into a trapezoid.

For a trapezoid with sides $\ell$ and $r$, and width $w$, the top side originally lies on the line $f_o(x)=\frac{r-\ell}{w}\cdot x + \ell$. To pinch the top-right corner to a new height $h$, it has to lie on the line $f_r(x)=\frac{h-\ell}{w}\cdot x + \ell$. Therefore, we define the required (non-linear) transformation as $Pinch_r(x,y) = \left(x,\frac{f_r(x)}{f_o(x)}\cdot y\right)$. A similar function $Pinch_\ell(x,y)$ can pinch a top-left vertex by using $f_\ell(x)=\frac{r-h}{w}\cdot x + h$.

To calculate the value of $h$, consider a composed diagram with top vertices at $V_1=(0,\ell_1)$ and $V_2=(w_1+w_2,r_2)$. The top-right vertex of D1 and top-left vertex of D2, both initially at $(w_1,r_1)=(w_1,\ell_2)$, must now become collinear with $V_1$ and $V_2$. That line's equation is $y=\frac{r_2-\ell_1}{w_1+w_2}\cdot x + \ell_1$. Thus, the new middle vertex at $x=w_1$ has a height of $y=\frac{w_1r_2+w_2\ell_1}{w_1+w_2}$, which is the desired $h$.

\subsection{Parallel composition (tensoring)}

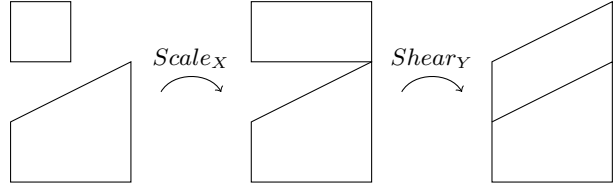
\begin{figure}[t]
    \centering
    \begin{tikzpicture}[scale=0.8]
        \draw (0,0) -- (0,1) -- (2,2) -- (2,0) -- cycle;
        \draw (0,2) -- (0,3) -- (1,3) -- (1,2) -- cycle;
        \draw[->] (2.5,1.5) to[out=60, in=120] node[midway,above] {$Scale_X$} (3.5,1.5);
    \begin{scope}[shift={(4,0)}]
        \draw (0,0) -- (0,1) -- (2,2) -- (2,0) -- cycle;
        \draw (0,2) -- (0,3) -- (2,3) -- (2,2) -- cycle;
        \draw[->] (2.5,1.5) to[out=60, in=120] node[midway,above] {$Shear_Y$} (3.5,1.5);
    \end{scope}
    \begin{scope}[shift={(8,0)}]
        \draw (0,0) -- (0,1) -- (2,2) -- (2,0) -- cycle;
        \draw (0,1) -- (0,2) -- (2,3) -- (2,2);
    \end{scope}
    \end{tikzpicture}
    \caption{Transformations to tensor two right trapezoids}
    \label{fig:tensoring-brickdiagrams}
    \end{figure}

When tensoring two right trapezoids (D1$\otimes$D2), the bottom side of D1 needs to be reshaped to match the slanted top side of D2. This can be done using a two-step (linear) transform on D1, as shown in \autoref{fig:tensoring-brickdiagrams}. The required transformation is the composite $Shear_Y\ \circ\ Scale_X$, where $ Scale_X(x,y) = \left(\frac{w_2}{w_1}\cdot x,y\right)$ scales $x$ by the ratio of the bottom to top widths, fixing $y$; and $ Shear_Y(x,y) = \left(x,\frac{\ell_2 - r_2}{w_2}\cdot x + y\right) $ shears $y$ by an amount proportional to the arity difference over the new width, fixing $x$.

For aesthetic and practical reasons, we prefer to scale both trapezoids to the maximum of widths $w_1$ and $w_2$ before shearing and tensoring. By scaling to the widest of the two, details are not unintentionally compressed or overlapped due to cramming into a smaller area.

\section{Implementation}

We implemented this method in the programming language Haskell, using a graphics library called \texttt{Diagrams} \cite{diagrams-manual}. The source code repository can be found at \url{https://github.com/celrm/stringdiagrams}.

The \texttt{Diagrams} library was a major inspiration for this work, but played a smaller role in the final implementation than originally expected because of a certain missing feature. \texttt{Diagrams} is a \emph{declarative} graphics library, in which the source code of a diagram is intended to roughly reflect its geometry. The library provides a datatype representing diagrams and \emph{combinators} for composing existing diagrams to build new diagrams. Among these are the operators \texttt{|||} and \texttt{===}, which join two diagrams horizontally and vertically respectively. An idealised implementation is depicted in \autoref{fig:naive-haskell}.

\begin{figure}[t]
\begin{verbatim}
render (Leaf l)           = renderLeaf l
render (Sequential d1 d2) = (render d1) ||| (render d2)
render (Parallel d1 d2)   = (render d1)
                                ===
                            (render d2)
\end{verbatim}
\caption{Idealised Haskell implementation}
\label{fig:naive-haskell}
\end{figure}

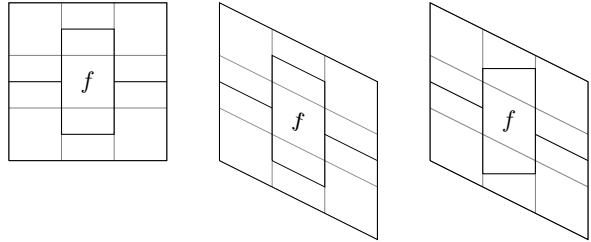
\begin{figure}[t]
\[ \begin{tikzpicture}[scale=0.7]
    \begin{scope}
	\draw[help lines] (0, 0) grid (3, 3);
	\draw (0, 0) to (0, 3) to (3, 3) to (3, 0) to (0, 0);
	\draw (1, .5) to (1, 2.5) to (2, 2.5) to (2, .5) to (1, .5);
	\node at (1.5, 1.5) {$f$};
	\draw (0, 1.5) to (1, 1.5); \draw (2, 1.5) to (3, 1.5);
\end{scope}
\begin{scope}[transform canvas={yslant=-.5}, xshift=4cm, yshift=2cm]
	\draw[help lines] (0, 0) grid (3, 3);
	\draw (0, 0) to (0, 3) to (3, 3) to (3, 0) to (0, 0);
	\draw (1, .5) to (1, 2.5) to (2, 2.5) to (2, .5) to (1, .5);
	\node at (1.5, 1.5) {$f$};
	\draw (0, 1.5) to (1, 1.5); \draw (2, 1.5) to (3, 1.5);
\end{scope} 
\begin{scope}[xshift=8cm]
	\draw[help lines, yslant=-.5] (0, 0) grid (3, 3);
	\draw[yslant=-.5] (0, 0) to (0, 3) to (3, 3) to (3, 0) to (0, 0);
	\draw (1, -.25) to (1, 1.75) to (2, 1.75) to (2, -.25) to (1, -.25);
	\node at (1.5, .75) {$f$};
	\draw (0, 1.5) to (1, 1); \draw (2, 0.5) to (3, 0);
\end{scope}
\end{tikzpicture} \vspace*{-1.5cm} \]
\vspace{1cm}
\caption{An example diagram (left) with a shear transformation applied to its contents (middle) and its coordinate system (right)}
\label{fig:coordinate-transform}
\end{figure}

Another feature of \texttt{Diagrams} that inspired us is the ability to apply transformations to a diagram's \emph{coordinate system} rather than its \emph{contents}. This is depicted in \autoref{fig:coordinate-transform}. Unfortunately this feature currently only works with very simple graphical elements, which does not include text or curves. Adding this feature to the library would simplify our implementation.

In the existing implementation, additional graphical elements like text labels or boxes must be managed separately from the diagram's coordinates, which are defined by an explicit trapezoidal bounding box; transformations have thus to be applied to these lists heterogeneously. However, wrapping these details in custom classes allows us to use a very similar syntax to the idealised implementation. 
The only requirement is then a method for converting our custom datatypes to the Diagram type from the library to be rendered. This out-of-the-box type cannot be defined as a direct instance of our custom class because it cannot be non-linearly deformed \textemdash which we need for horizontal composition. This non-homogenous approach also brings some advantages, such as the ability to smooth out horizontal connections at a flat angle, for a more aesthetically pleasing result.

The final diagram is generated by a simple recursive function without knowing the general placement of any graphical element beforehand. This is a major advantage of our method, compared to others which rely on pre-calculated graph layouts. We can also use this heterogeneity to easily carry semantics along with our diagrams, which get compiled at the same time as their graphic representation is rendered. In fact, this is a practical proof that drawing is merely another form of compilation. We illustrate this by additionally supporting a matrix algebra backend, where each diagram denotes a matrix, composition is matrix multiplication, and tensoring is the direct sum of matrices.

\begin{figure}[t]
\centering
\includegraphics[width=0.9\linewidth]{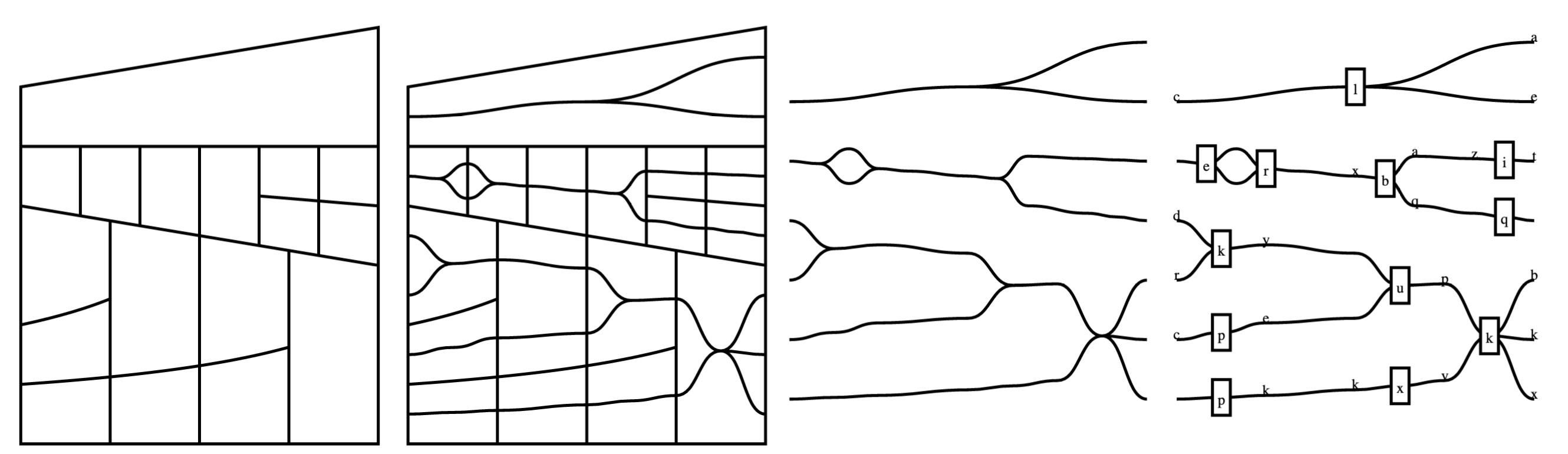}
\caption{Randomly-generated test example}
\label{fig:test-example}
\end{figure}

An example of the output generated by our code on a randomly-generated test example is depicted in \autoref{fig:test-example}.
The most similar existing software to ours is the rendering part of DisCoPy \cite{discopy}, which has similar capabilities but uses a different data representation and does not use recursive rendering.

\bibliographystyle{splncs04}
\bibliography{rendering-diagrams}

\end{document}